\documentclass[11pt,a4paper,reqno]{amsart} %maths

\usepackage{amssymb,graphicx}%maths
\usepackage{amssymb}

\newcommand{\koniec}{\begin{flushright}  $\Box $ \end{flushright}}
\newtheorem{theo}{Theorem}[section] 
\newtheorem{prop}[theo]{Proposition}  
\newtheorem{lemma}[theo]{Lemma}
\newtheorem{defi}[theo]{Definition}

\def\theequation{\thesection.\arabic{equation}}

\newcounter{mnotecount}[section]

\renewcommand{\themnotecount}{\thesection.\arabic{mnotecount}}

\newcommand{\mnote}[1]%{}%
{\protect{\stepcounter{mnotecount}}$^{\mbox{\footnotesize
$%\!\!\!\!\!\!\,
\bullet$\themnotecount}}$ \marginpar{%\color{red}%
\raggedright\tiny\em
$\!\!\!\!\!\!\,\bullet$\themnotecount: #1} }

\newcommand{\hook}{{\setlength{\unitlength}{11pt}   % adjust pt size here
                   \begin{picture}(.833,.8)
                   \put(.15,.08){\line(1,0){.35}}
                   \put(.5,.08){\line(0,1){.5}}
                   \end{picture}}}
\newcommand{\CP}{\mathbb{CP}}
\newcommand{\C}{\mathbb{C}}

\newcommand{\R}{\mathbb{R}}
\newcommand{\Z}{\mathbb{Z}}

\def\p{\partial}
\def\be{\begin{equation}}
\def\ee{\end{equation}}

\def\bea{\begin{eqnarray}}
\def\eea{\end{eqnarray}}

\begin{document}\date{April 25, 2016}
%%%%%%%%%%%%%%%%%%%%%%%%%%%%%%%%%%%%%%%%
\title{On the quadratic invariant of binary sextics}
%%%%%%%%%%%%%%%%%%%%%%%%%%%%%%%%%%%%%%%%
\author{}
\author{Maciej Dunajski}
\address{Department of Applied Mathematics and Theoretical Physics\\ 
University of Cambridge\\ Wilberforce Road, Cambridge CB3 0WA\\ UK.}
\email{m.dunajski@damtp.cam.ac.uk}
\author{Roger Penrose}
\address{The Mathematical Institute\\
University of Oxford,\\
Andrew Wiles Building,\\
Woodstock Road, Oxford OX2 6GG
\\ UK.
}
\email{penroad@wadh.ox.ac.uk}

%%%%%%%%%%%%%%%%%%%%%%%%%%%%%%%%%%%%%%%%%
\begin{abstract} 
We provide a geometric characterisation of binary sextics
with vanishing quadratic invariant. 
%%%%%%%%%%%%%%%%%%%%%%%%%%%%%%%%%%%%%%%%%%
\end{abstract}   
\maketitle
%%%%%%%%%%%%%%%%%%%%%%%%%%%%%%%%%%%%%%%%% 
%%%%%%%%%%%%%%%%%%%%%%%%%%%%%%%%%%%%%%%%%%
\section{Introduction}
Classical invariant theory was formulated and developed
by Cayley, Salmon, Sylvester, Young and others in the second half
of the $19^{th}$ century. After a quiet period which lasted for the most
part of the $20^{th}$ century, the theory has reappeared in
algebraic geometry as the {\em geometric invariant theory} \cite{mum}, in representation theory  of $GL(2)$, and in other
areas of modern mathematics \cite{dolg}. Each of these areas replaces
the classical terminology  by its own language.
In this paper we shall nevertheless follow the old--fashioned terminology
of \cite{Grace_Young}. This is in line with other modern expositions
of the subject \cite{olver,kung}. 

 Despite the recent developments, some of the problems left over
from the $19^{th}$ century remain open. One class of such problems
has to do with finding the interpretation of the vanishing of  invariants and covariants. This should be expressed in terms of the underlying projective geometry of roots on the Riemann sphere.
The aim of this paper is to solve a problem from this category, and
interpret the vanishing of the {\it apolar} invariant for a binary sextic, and more generally for  binary quantics of even degree.

Consider a polynomial of degree six in $x$ with complex
coefficients ${\boldsymbol{\psi}}=(\psi_0, \dots, \psi_6)$
\be
 \label{sextic}
\psi(x)=\psi_0x^6+6\psi_1x^5+15\psi_2x^4+20\psi_3x^3+15\psi_4x^2+
6\psi_5x+\psi_6.
\ee
Substituting
\[
x=\frac{a \tilde{x}+b}{c\tilde{x}+d},\quad
\mbox{where}\quad ad-bc=1,
\]
and multiplying the resulting expression by 
$({c\tilde{x}+d})^6$ to clear the denominators gives
the polynomial $\tilde{\psi}(\tilde{x})$
%\[
%\tilde{\psi}_0x^6+6\tilde{\psi}_1x^5+15\tilde{\psi}_2x^4+20\tilde{\psi}%_3x^3
%+15\tilde{\psi}_4x^2+6\tilde{\psi}_5x+\tilde{\psi}_6,
%\]
with the coefficients ${\boldsymbol{\tilde\psi}}$ given by a linear 
transformation of ${\boldsymbol{\psi}}$. The function
\be
\label{inv1}
{\mathcal{I}}({\boldsymbol{\psi}})=2\psi_0\psi_6-12\psi_1\psi_5+30\psi_2\psi_4-20\psi_3^2
\ee
is an {\em invariant} of the sextic, as  ${\mathcal{I}}(\boldsymbol{\psi})={\mathcal{I}}(\boldsymbol{\tilde{\psi}})$. 
 
 Vanishing of any invariant of a binary quantic -- the precise definitions are given in the next Section -- describes some geometric property of the configurations of the roots of the quantic
regarded as points on the two--dimensional sphere $\CP^1=\C+\{\infty\}$. The analog of the quadratic invariant (\ref{inv1}) 
can be constructed for any polynomial of even degree - see formula 
(\ref{qin}). 
In this paper we find a geometric interpretation of the condition ${\mathcal I}=0$.
This natural question does not seem to have been answered by the classical invariant theorists in the $19^{th}$ and early $20^{th}$ centuries, except when the the quantic has degree two or four: A generic quartic has four distinct roots, and the  condition $\mathcal{I}=0$ implies that their cross  ratio  is a cube root of unity. This is the equianharmonic condition. The roots of the quartic, when viewed as points
on the Riemann sphere, can in this case be 
transformed into vertices of a regular
tetrahedron by a M\"obius transformation. The answer for binary sextics
appears not to be known, and trying to understand this case in the context of twistor theory of $G_2$ structures
\cite{DG10, DD} is one motivation for this paper. It turns out that the sextic
case can be reduced to the quartic in a sense made precise by the following Theorem
\begin{theo}
\label{main_thm}
Let $X_1, X_2, X_3, X_4$ be four points on a two--dimensional sphere such that the stereographic projection of one of the roots of the sextic {\em(\ref{sextic})}
from any of these four points lies in the centroid of the projections
of the remaining five roots. Then ${\mathcal{I}}({\boldsymbol{\psi}})=0$ if and only if the points $X_1, \dots, X_4$
can be transformed into vertices of a regular tetrahedron by a M\"obius
transformation (if they are distinct), or if at least three
of these points coincide.
\end{theo}
First we will need to establish that there are indeed only four points,
up to multiplicity, with the property stated in this Theorem. 
This, together with the rest of the proof with be presented in Section
\ref{section4}. 

 The invariant (\ref{inv1}) is quadratic in the coefficients
of the sextic, and thus given any five points on the sphere there will generically exist two points which complement these five points
to roots of a sextic with $\mathcal{I}=0$. In Proposition \ref{non_ger_prop}
we shall characterise the non--generic configurations of five points
such that any choice of a distinct sixth point  yields a sextic
with $\mathcal{I}\neq 0$. It will be shown that any such
non-generic configuration is projectively equivalent to vertices
of a square pyramid.

 The sextic case  is special in some ways, but
the general method in the paper together with an inductive argument 
applies to binary quantics of any even degree.
\subsubsection*{Acknowledgements}
We thank Robert Bryant, Mike Eastwood, Nigel Hitchin and others 
for useful discussions.
%%%%%%%%%%%%%%%%%%%%%%%%%%%%%%%%%%%%%%%%%%%%%%%%%%%%%%%%%%%%%%%%%%%%%%
\section{Quantics and invariants}
\label{section2}
A binary quantic is a homogeneous polynomial in two variables which we
shall call $(x, y)$. We shall consider binary quantics of even degree
\be
\label{quantic1}
\psi(x, y)=\sum_{k=0}^{2n}{2n\choose k} \psi_k x^{2n-k}y^k.
\ee
The coefficients of the quantic ${\boldsymbol{\psi}}=(\psi_0, \dots, \psi_{2n})$ are assumed to be complex numbers. There exists a unique, up to an overall scale, quadratic invariant
\be
\label{qin}
\mathcal{I}(\psi)=2\sum_{k=0}^{2n} (-1)^{2n-k} {2n\choose k}\psi_k\psi_{2n-k}
\ee
one can associate to 
the quantic\footnote{Sylvester \cite{sylvester}
calls it the {\em quadrinvariant}.
In his latter works, see e.g. \cite{atom}, he proposed an analogy between classical
invariant theory and molecular chemistry. A binary sextic would
correspond to an atom with six free valent electrons, and the 
quadratic invariant is then the bi-atomic molecule.}.
The {invariance} of $\mathcal{I}$ is to be understood in the following way: Consider the linear action of $GL(2, \C)$ on $\C^2$ given 
by the change of variables
\[
{x}=a \tilde{x}+b \tilde{y}, \qquad {y}=c \tilde{x}+ d \tilde{y}, \qquad
ad-bc\neq 0.
\]
Given  a binary quantic $\psi(x, y)$, 
let $\widetilde{\psi}({\tilde x}, {\tilde y})$ be a binary quantic 
given by
\begin{eqnarray*}
\widetilde{\psi}({\tilde x}, {\tilde y})&=&
\sum_{k=0}^{2n}{2n\choose k}\;\psi_k\;
  (a \tilde{x}+b \tilde{y})^{2n-k}\;(c \tilde{x}+d \tilde{y})^k\\
&=&\tilde{\psi}_0\tilde{x}^{2n}+2n\tilde{\psi}_1 \tilde{x}^{2n-1}\tilde{y}+n(2n-1) \tilde{\psi}_2 \tilde{x}^{2n-2}\tilde{y}^2+
\dots +\tilde{\psi}_{2n} \tilde{y}^{2n}.
\end{eqnarray*}
This induces an irreducible embedding  $GL(2, \C)\subset GL(2n+1, \C)$, as
the $(2n+1)$ coefficients of $\tilde{\psi}$
are linear homogeneous functions of the coefficients of $\psi$.
\begin{defi}
A {covariant} of a binary quantic 
is a polynomial $I=I(\psi_0, \dots, \psi_{2n}, x, y)$ 
such that
\[
I( \tilde{\psi}_0, \dots, \tilde{\psi}_{2n}, \tilde{x}, \tilde{y})
=(ad-bc)^w I( \psi_0, \dots, \psi_{2n}, x, y ).
\]
The number $w$ is called the weight of the covariant.
A covariant which only depends on the coefficients of the quantic, and
not on $(x, y)$ is called an invariant. 
\end{defi}\noindent
Thus the degree--two  invariant (\ref{qin})
has weight $2n$. There are other invariants of degree higher than two \cite{Grace_Young}.
In the case of the binary 
sextic
there are four more invariants, of degree 4, 6, 10 and 15 respectively 
connected with $\mathcal{I}$ by a syzygy of degree 30.
\subsection{Transvectants}
\label{section3}
Let $V_{m}=\mbox{Sym}^{m}(\C^{2\vee})$ be the $(m+1)$--dimensional complex vector space of binary quantics of degree $m$, where $\C^{2\vee}$ is the dual of $\C^2$.
Given two binary quantics $\phi\in V_n$ and $\psi\in V_m$, the $k$th transvectant is a map $<\;,\;>_k:V_m\times V_n\rightarrow
V_{m+n-2k}$ given by a quantic of degree $n+m-2k$
\be
\label{k_th_trans}
<\phi, \psi>_k:=\sum_{j=1}^k(-1)^j{k\choose j}\frac{\p^k \phi}{\p x^{k-j}\p y^j}
\frac{\p^k \psi}{\p x^{j}\p y^{k-j}}.
%\phi_{A_1\dots A_k (BC\dots D}{\psi^{A_1\dots A_k}}_{EF\dots G)}
\ee
%\[
%<\psi, \phi>_k=\varepsilon^{A_1B_1}\dots \varepsilon^{A_kB_k}
%\psi_{A_1\dots A_kA_{k+1}\dots A_m}\phi_{B_1\dots B_kB_{k+1}\dots B_n}\pi^{A_{k+1}}
%\dots \pi^{A_m}\pi^{B_{k+1}}\dots\pi^{B_n}.
%\]
Thus, for any $k\leq\mbox{min}(m, n)$,  transvectants are covariants of weight $k$ and degree two.
One of the results in the classical invariant theory is that
all covariants and invariants arise from the transvectants 
operations \cite{Grace_Young,olver}.
\begin{defi}
\label{defi_apolar}
The quantic $\phi\in V_n$ is 
{\it apolar} to $\psi\in V_m$ where $m\geq n$ if $<\psi, \phi>_n=0$. 
\end{defi}\noindent
%Thus, in particular,
%the quadratic invariant (\ref{qin}) vanishes if the quantic $\psi$ is %apolar to itself.
% The apolar
%covariant of $\psi$ and $\phi$ is the $n$th transvectant
%\[
%<\psi, \phi>_n=\psi_{A_1A_2\dots A_m}\phi^{A_1A_2\dots A_n} \;
%\pi^{A_{n+1}}\dots \pi^{A_{m}}.
%\]
In the special case $n=m$ the apolarity condition is given by 
\[
\sum_{k=0}^n (-1)^{n-k} {n\choose k}\psi_k\phi_{n-k}=0.
\]
Any quantic of an odd degree is apolar to itself.
A quantic of an even degree is apolar to itself iff the quadratic invariant $\mathcal{I}$ given by  (\ref{qin})
vanishes.
%%%%%%%%%%%%%%%%%%%%%%%%%%%%%%%%%%%%%%%%%%%%%%%%%%%%%%%%%%
\section{Characterisation of $\mathcal{I}$}
\label{section4}
Let $\psi$  be a sextic (\ref{sextic}) which is generic in the sense that its six complex roots are distinct,
and let ${\mathcal{I}(\psi)}$ be given by (\ref{inv1}).
\begin{prop}
\label{prop11}
Let $P_1, \dots, P_6$ be six points on the sphere $\CP^1$ corresponding to
the roots of a sextic $\psi$. Then $\mathcal{I}(\psi)=0$ if and only if the four roots
of the quartic $<\kappa, \rho>_1$ can be transformed into vertices
of a regular tetrahedron (if they are distinct) or contain a root of multiplicity at least three. Here $\kappa$ and $\rho$ are any quintic
and a linear form respectively such that $\psi(x)=\kappa(x)\rho(x)$.
\end{prop}\noindent
{\bf Proof.}  Any quintic corresponding to distinct points 
$P_1, \dots, P_5$ associates
four points on the sphere to a given point $P_6$, such that
the quartic defining the four points
is a transvectant of the quintic with a linear form corresponding to $P_6$.
Given $\psi(x, y)=\kappa(x, y)\rho(x, y) $, and using $<\rho, \rho>_1=0$ we compute 
\[
\mathcal{I}(\psi)=<\psi, \psi>_6=-<\delta, \delta>_4=-\mathcal{I}(\delta), \quad \mbox{where}\quad \delta=<\kappa, \rho>_1.
\]
Thus ${\mathcal{I}}(\psi)=0$ iff ${\mathcal{I}}(\delta)=0$, where
$\delta$ is the quartic defined above.  But ${\mathcal{I}}(\delta)=0$
if at least three of the roots of $\delta$ coincide, or if the
cross ratio of the four distinct roots is a cube root of unity
(the equianharmonic case). For completeness we give the proof
of this fact, well known in the $19^{th}$ century literature.
Set $y=1$, consider a general quartic
\[
\delta=\delta_0 x^4+4\delta_1 x^3+ 6\delta_2 x^2+ 4\delta_3 x +\delta_4
\]
with four distinct roots. This has 
$\mathcal{I}(\delta)=2\delta_0\delta_4 -8\delta_1\delta_3+6(\delta_2)^2$.
Now consider the quartic corresponding to a regular tetrahedron with one
vertex at $\infty$. This quartic is represented by $(x-1)(x-\omega)(x-\omega^2)$
where $\omega^3=1$, and we find $\mathcal{I}=0$. The cross ratio
of the roots is given by $\omega$, and conversely any quartic
with the cross ratio of four roots given by a cube root of unity
is projectively equivalent to the tetrahedral quartic.
\koniec
We can now give the proof of the result stated in the Introduction.
In the proof, and in the remaining part of the paper we
shall write
\be
\label{notation}
\psi=P_1P_2\dots P_m
\ee
to denote a binary quantic defined (up to an overall non--zero multiple)
by its roots corresponding to points $P_1, P_2, \dots, P_m$ on the sphere.\\
{\bf Proof of Theorem \ref{main_thm}.} 
Consider the quartic $\delta=<\kappa, \rho>_1$ introduced in Proposition
\ref{prop11}. Let its four roots correspond to the (not necessarily distinct) points $X_1, \dots, X_4$. 
Thus, using the notation (\ref{notation}), the quartic $\delta$ is apolar
(in the sense of Definition \ref{defi_apolar})
to a quartic $XXXX$, where $X$ is any of the roots $X_i$.
\begin{center}
\includegraphics[width=8cm,height=4cm,angle=0]{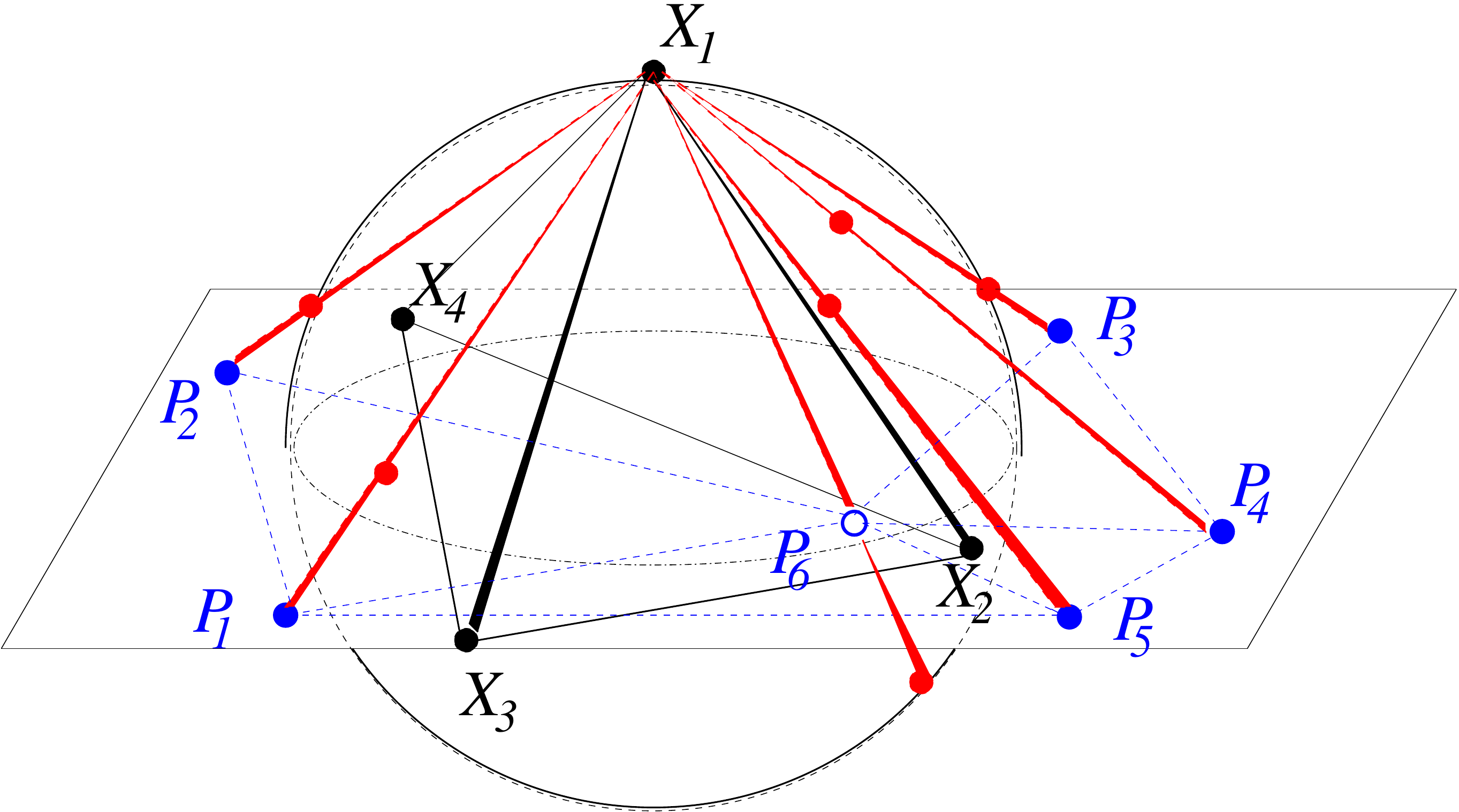}
\begin{center}
{{\bf Figure 1.} {\em Stereographic projection of $P_6$ to the centroid
of $P_1, \dots, P_5$.}}
\end{center}
\end{center}
Equivalently
the quintic $\kappa=P_1P_2P_3P_4P_5$ is apolar to the quintic
$\chi=P_6X_iX_iX_iX_i$ for any $i=1, \dots, 4$. This happens if and
only if the stereographic projection of $P_6$ from any of the four 
points $X_i$ lies in the centroid of the projections of the five
points $P_1, \dots, P_5$. This fact was observed in \cite{penrose, ZP94}.
Indeed,
if the homogeneous coordinates of the projected points $P_i$ 
are $(1, -x_i)$ where $i=1, \dots, 6$ and the coordinates of the north pole $X$  are $(0, 1)$ then 
$<P_i, X>=1, <P_i, P_j>=x_i-x_j$, so
\[
<\kappa, \chi>_5= \sum_{i=1}^5(x_i-x_6)=0, \qquad\mbox{and}\quad x_6=\frac{1}{5}\sum_{i=1}^5 x_i.
\]
We note that the centroid of the set of points is not invariant under the projective transformations, but is invariant under the subgroup of affine
transformations $x\rightarrow \alpha x+\beta$ preserving the north pole of the stereographic projection. The statement in the Theorem \ref{main_thm}
is nevertheless projectively invariant.
\koniec
This proof extends to binary quantics of arbitrary even degree. Assume we know a geometrical interpretation for vanishing of $\mathcal{I}$ for a binary quantic $\psi_{2n-2}$ of degree $2n-2$. We can now characterise the self--apolarity of a binary quantic $\psi_{2n}$ of degree $2n$ by the following inductive argument:  Given $2n-1$ distinct points $(P_1, \dots, P_{2n-1})$ on the sphere, let $\kappa$
be the unique (up to a non--zero multiple) binary quantic 
of degree $(2n-1)$
with these points corresponding to its roots. The points 
$(P_1, \dots, P_{2n-1}, P_{2n})$ are roots of a binary quantic $\psi_{2n}$ with ${\mathcal I}(\psi_{2n})=0$
if a linear form $\rho$ with a root corresponding to the point $P_{2n}$ is such that
the binary quantic $\psi_{2n-2}=<\kappa, \rho>$ has $\mathcal{I}(\psi_{2n-2})=0$.

 Given $2n$ unordered points $\psi_{2n}=\{P_1, P_2, \dots, P_{2n}\}$ on the sphere, 
split them into 
a set of $2n-1$ points $\kappa=\{P_1, P_2, \dots, P_{2n-1}\}$
together with one point $\{P_{2n}\}$. 
Let $\delta=\{X_1, X_2, \dots, X_{2n-2}\}$ be the points on the sphere
such that the stereographic projection of $P_{2n}$ from any $X_i$
is the centroid of the stereographic projections
of $P_1, \dots, P_{2n-1}$. Then the set of points $\psi$ is self-apolar
iff the set of points $\delta$ is self--apolar.
The set $\delta$ consist of the roots of the polynomial equation
\be
\label{polyn}
%\kappa_{ABC\dots D}\rho^A X^BX^C\dots X^D=0
<\kappa, \chi>_{2n-1}=0,\quad\mbox{where},\quad \chi=P_{2n}XX\dots X
\ee
counted with multiplicity.
\subsubsection*{Example.} 
The multiplicities of the elements of $\delta$ can depend on the choice of the point $P_{2n}$ from the set $\psi$.
Consider the sextic 
\be
\label{sextic_5}
\psi=(x-1)(x-\omega)(x-\omega^2)(x-\omega^3)(x-\omega^4), \quad\mbox{where}\quad \omega^5=1
\ee
corresponding to a pentagonal pyramid with one root placed at $\infty$. This has ${\mathcal I}=0$.
Thus $\psi=\{1, \omega, \omega^2, \omega^3, \omega^4, \infty\}$. Taking $P_6=\infty$
gives the quartic (\ref{polyn}) to be $x^4$ which has one quadrupole 
root $x=0$ and thus is self--apolar. Choosing instead $P_6=1$ gives the polynomial (\ref{polyn}) 
\[
x^4+6x^3+6x^2+6x+6=0
\]
which has four distinct roots with equianharmonic cross ratio.
\subsection{Canonical form}
Almost all sextics can be put in the Sylvester 
Canonical Form \cite{sylvester}
\be
\label{can_form_syl}
\psi=C u^6 + A v^6 + B w^6 + uvw(u-v)(v-w)(w-u),
\ee
where $u+w+v=0$ - thus to obtain (\ref{sextic}) set $u=x, v=1, w=-x-1$ in this formula\footnote{The exceptional sextics which can not be put in
the canonical form are classified in \cite{Eastwood}.}.
In this form the quadratic invariant (\ref{inv1})
is 
\[
\mathcal{I}({\boldsymbol{\psi}})=2CA + 2CB + 2BA - 2.
\]
If the sextic does not have a root at $\infty$ we can assume that
$C+B\neq 0$, and solve for $A=(1-CB)/(C+B)$. This gives a canonical
form of a generic self-apolar sextic. It
parametrises a non--singular open orbit in the space of all self-apolar sextics.
Setting $B=(b+c)/(3b-3c), C=(6-b-c)/(3b-3c)$ gives 
\be
\label{sylv_can}
\psi=x^6+2bx^5+5bx^4+\frac{20}{6}(b+c)x^3+5cx^2+2cx+
\frac{1}{36}(b+c)^2+\frac{1}{4}(b-c)^2,
\ee
so the self-apolar sextic depends on two parameters.
In general  the binary quantic of degree $2n$ with
$\mathcal{I}=0$ depends on $2n-4$ arbitrary parameters, up to the M\"obius transformation.
\subsection{Catalectant}
We shall finish off this section giving an alternative interpretation
of the condition ${\mathcal I}=0$, which brings up another quartic invariant of binary sextics. 

The seven--dimensional complex vector space of binary sextics $V_{6}$ belongs to the space
of endomorphisms of the four--dimensional  space of binary cubics $V_3$, where the endomorphism
corresponding to a sextic $\psi$ is given by the
transvectant
$
\phi\rightarrow <\psi, \phi>_3,
$
where $\phi\in V_3$.
%The resulting endomorphisms are 
%traceless. The space of traceless endomorphisms of $V_n$ has dimension
%$({(n+1)(n+2)})/{2}-1$,
%which is greater than the dimension of $V_{2n}$  unless $n=2$.
%\vskip5pt
%where we expect two conditions 
%characterising  a general traceless endomorphism which arises from a sextic. 
Consider a complex eigenvalue problem
\be
\label{map}
<\psi, \phi>_3=\lambda\phi.
%{\psi_{ABC}}^{DEF}\phi_{DEF}=\lambda\phi_{ABC}.
\ee
The eigenvalue $\lambda$ has weight three under the $GL(2)$ action on $V_6$.
It is also an eigenvalue, in the ordinary sense, of a 4 by 4 
matrix
corresponding to $\psi$ by choosing a basis of $V_3$. 
If $\lambda_1, \lambda_2, \lambda_3, \lambda_4$ are the four eigenvalues
of $\psi$ then the characteristic polynomial is
%\begin{eqnarray*}
%\sum_{i}\lambda_i&=&{\psi_{ABC}}^{ABC}=0,\\
%\sum_{i}(\lambda_i)^2&=&{\psi_{ABC}}^{DEF}{\psi_{DEF}}^{ABC}=-\mathcal{I}%,\\
%\sum_{i}(\lambda_i)^3&=&{\psi_{ABC}}^{DEF}{\psi_{DEF}}^{GHI}{\psi_{GHI}}%^{ABC}
%=0,\\
%\sum_{i}(\lambda_i)^4&=&{\psi_{ABC}}^{DEF}{\psi_{DEF}}^{GHI}{\psi_{GHI}}%^{JKL}
%{\psi_{JKL}}^{ABC}=\mathcal{J}.
%\end{eqnarray*}
a quartic\footnote{This is the equation {\bf R} on page 24 in \cite{sylvester}, used to construct the 
canonical form (\ref{can_form_syl}).} 
%(the factor 8 is inserted for convenience)
\begin{eqnarray}
\label{quartic}
\chi_\psi(\lambda)&=&
8(\lambda-\lambda_1)(\lambda-\lambda_2)(\lambda-\lambda_3)(\lambda-\lambda_4)\nonumber\\
&=&8\lambda^4+4\mathcal{I}\;\lambda^2-\mathcal{J}.
\end{eqnarray}
Such quartic is canonically associated to every sextic, and has weight 
$12$. The degree four invariant 
$\mathcal{J}=<<\psi, \psi>_3, <\psi, \psi>_3>_6$  appearing in (\ref{quartic}) is the {\it catalectant} of the sextic. Its zero set
is the closure of the
locus of sextics expressible as the sum of three $6$th powers \cite{Elliott}.
Equivalently, the catalectant vanishes
if the sextic admits an apolar cubic \cite{Grace_Young}.
The latter result follows directly from  
setting  $\lambda=0$ in (\ref{quartic}). If $\mathcal{I}=0$ and $\mathcal{J}\neq 0$, then the roots of (\ref{quartic})  form a {\em harmonic} set. This 
corresponds to a square on an equator in $\CP^1$.
%%%%%%%%%%%%%%%%%%%%%%%%%%%%%%%%%%%%%%%%%%%%%%%%%%%%%%%%%%%%%%%
\section{Maximally separated quintics}
\label{section41}
In this section we shall consider a problem of recovering a 
sextic $P_1P_2\dots P_6$
with ${\mathcal I}=0$ from a quintic. Let us project stereographically the sphere to the complex plane
from one of the roots - say $P_6=\infty$.
Given four points
$P_1, \dots,  P_4$ on the plane, we can now look for a fifth point $P_5$ such that
$\mathcal{I}=0$. Rewriting the sextic (\ref{sextic})  with the root 
corresponding to $P_6$ at 
$\infty$ 
as
\[
\psi=(x-x_1)(x-x_2)(x-x_3)(x-x_4)(x-x_5)
\]
and comparing the coefficients of various powers of $x$ we find that
\be
\label{quadratic}
\mathcal{I}=a\;{x_5}^2+2b\;x_5+c,
\ee
where $(a, b, c)$ are polynomials in 
$(x_1, \dots, x_4)$. This expression has at most two roots $x_5\neq \infty$, 
thus
given five distinct points on $\CP^1$ there exist at most two points such that
the invariant $\mathcal{I}$ of the sextic defining the resulting six
distinct points vanishes. Generically, if $a\neq 0$, there will be two such 
points. 
\subsubsection*{Example} Consider
a quintic corresponding to four points on the base of an equatorial regular pentagon, and a point at infinity.
Equation (\ref{quadratic}) gives
%a pyramid with a regular pentagon  as the base.
%We find that $\mathcal{I}=0$ for this configuration. 
%This follows from writing the sextic with a root at infinity as
%\[
%\psi=(x-1)(x-\omega)(x-\omega^2)(x-\omega^3)(x-\omega^4)
%\]
%in which case
%\[
%\mathcal{I}=
%-\frac{1}{120}{\omega}^{2} \left( {\omega}^{4}+{\omega}^{3}+{\omega}^{2}+\omega+1 \right)  \left( 3\,{\omega}^{4}
%+7\,{\omega}^{3}+11\,{\omega}^{2}+7\,\omega+3 \right)  \left( \omega-1 \right) ^{4}.
%\]
%Thus $\mathcal{I}$ vanishes if $\omega^5=1$, but it also vanishes
%for two other choices of $\omega$ and their complex conjugates (these
%points do not lie on the equatorial circle). Assume $\omega^5=1$. 
two possibilities for the sixth s.t. ${\mathcal I}=0$,
one of which is gives rise to a pentagonal pyramid (\ref{sextic_5})(Figure 2).
%Taking
%\be
%\label{pent}
%\psi=(x-x_5)(x-\omega)(x-\omega^2)(x-\omega^3)(x-\omega^4), \quad
%\mbox{where}\quad \omega^5=1
%\ee
%we find a quadratic equation (\ref{quadratic}) for $x_5$ with roots $x_5=1$ and $x_5=-3$. 
\begin{center}
\includegraphics[width=8cm,height=4cm,angle=0]{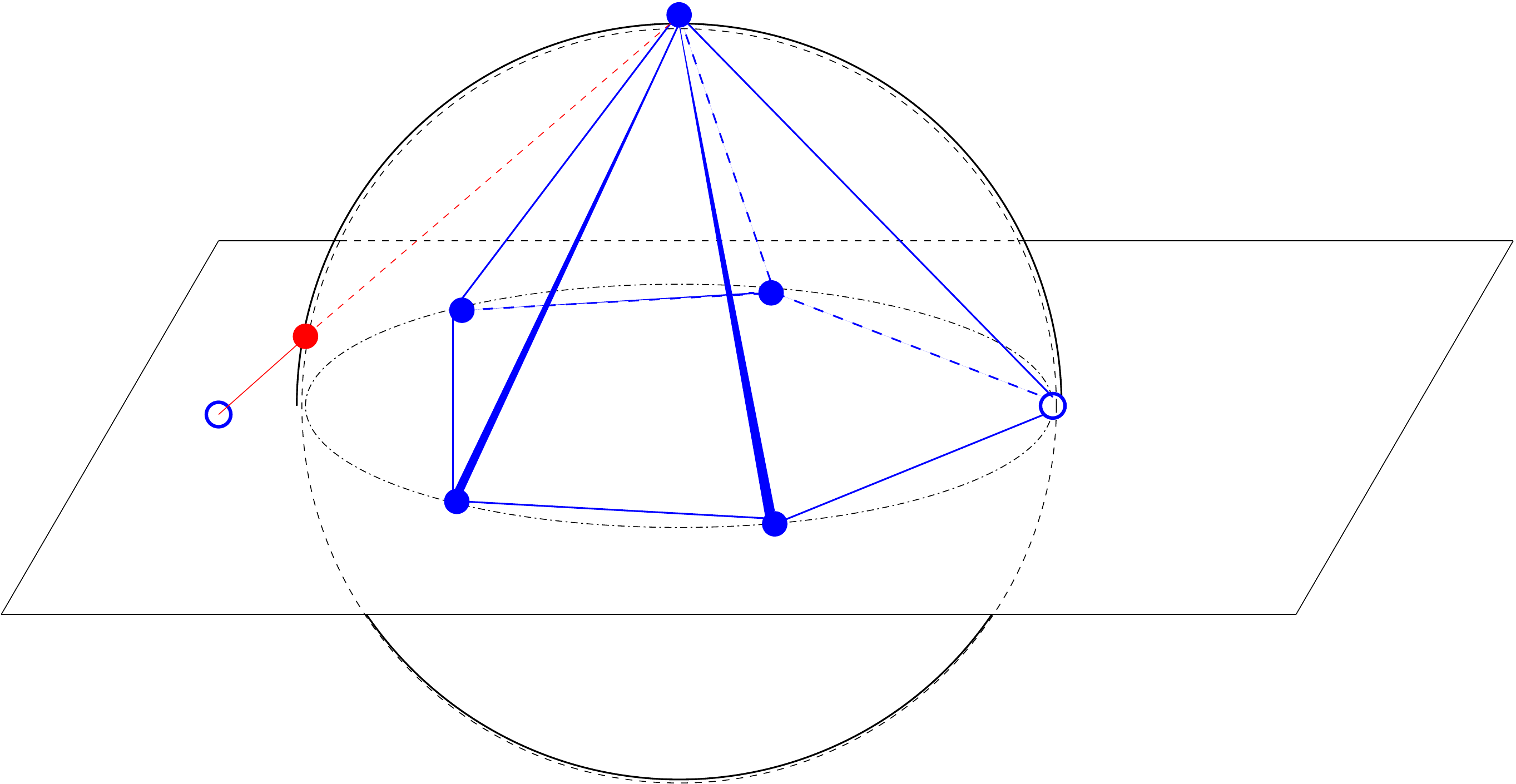}
\begin{center}
{{\bf Figure 2}. {\em Regular pentagon with $x_5=1$ and its companion with $x_5=-3$.}}
\end{center}
\end{center}
The next example is of a non-generic type
\subsubsection*{Example--Regular Octahedron.} Assume that the six points $P_i$ form the vertices
of the regular octahedron, and project from one of its vertices
$P_6$ which does not belong to the square. The resulting sextic
is $\psi=x^5-x$, and  $\mathcal{I}=1/3$. We can transform the four corners
of the square $P_1, \dots, P_4$ to $\pm 1, \pm i$, and
look for a point $P_5$ such that $\mathcal{I}=0$. The resulting sextic is
\[
(x^4-1)(x-x_5)
\]
and we find that $\mathcal{I}=1/3$ for any value of $x_5$. Thus
there is no $P_5$ which complements the five vertices
of the square pyramid to an octahedron with $\mathcal{I}=0$.
This example corresponds to $a=b=0, c\neq 0$ in (\ref{quadratic}).
\begin{center}
\includegraphics[width=4cm,height=8cm,angle=270]{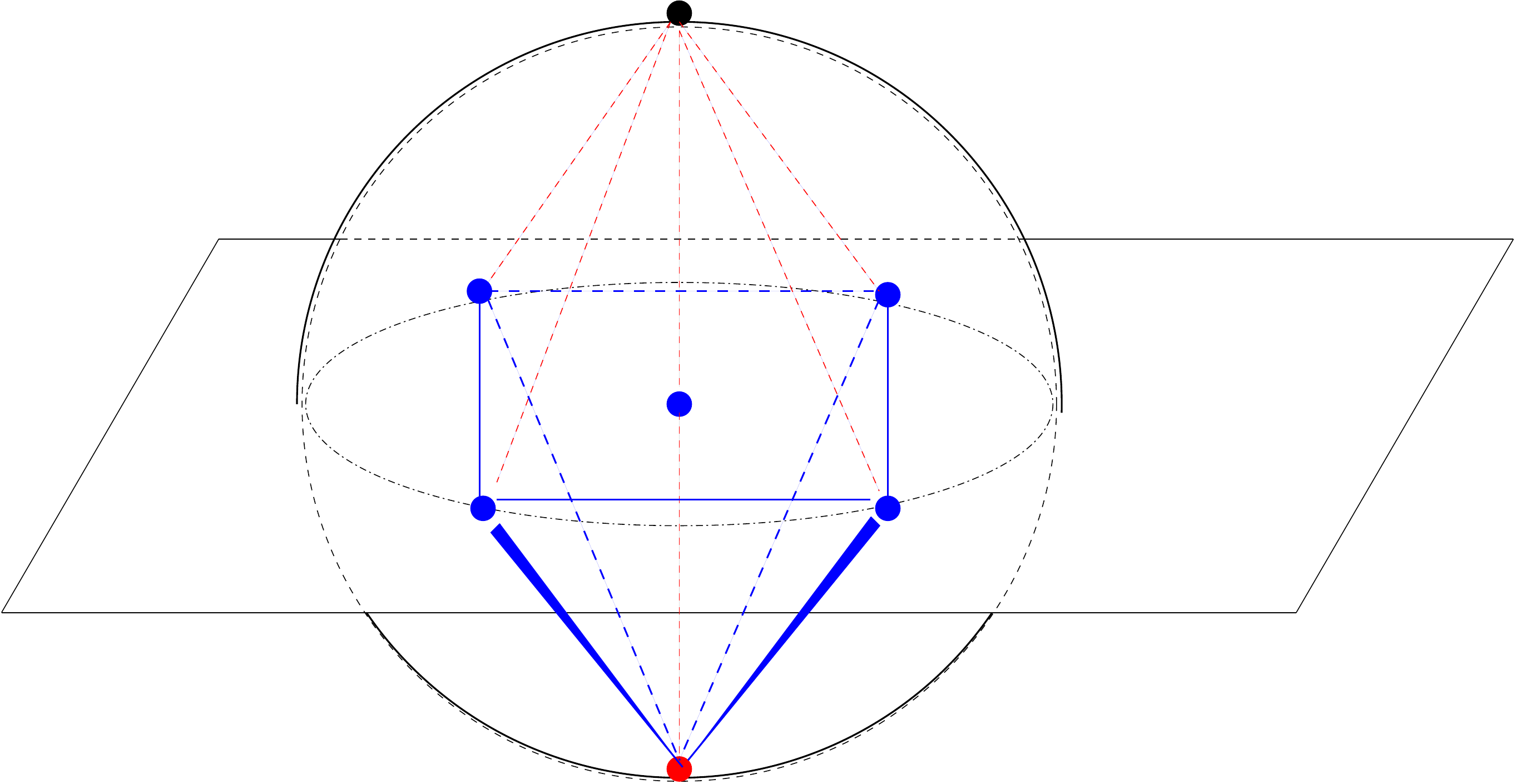}
\begin{center}
{{\bf Figure 3.} {\em Five maximally separated points.}}
\end{center}
\end{center}
\begin{defi}
We shall call five distinct points on the sphere {maximally separated}
if ${\mathcal{I}}\neq 0$ for the sextic with six distinct roots 
resulting from any choice of the sixth point.
\end{defi}
\begin{prop}
\label{non_ger_prop}
Any maximally separated five points on the sphere
can be transformed into vertices of a square pyramid by a M\"obius
transformation.
\end{prop}
{\bf Proof.}
Stereographicaly project from any 
of the five  maximally separated points. The remaining four points
are the distinct roots of a quartic 
\[
\gamma(x)=x^4+4b_1x^3+6b_2x^2+4b_3x+b_4.
\]
Now consider the sextic $\psi$ with one root at infinity given by
$\psi=\gamma(x)(x-x_5)$. Computing (\ref{quadratic}) and equating 
$a$ and $b$ to zero gives $b_2=b_1^2, b_3=b_1^3$, and
now $\mathcal{I}\neq 0$ unless all four roots of $\gamma$ coincide which 
we have excluded. The roots of the resulting quartic
are of the form $\alpha\pm\beta, \alpha\pm i\beta$, where
$\alpha$ and $\beta$ are complex numbers depending on $b_1$ and $b_4$.
These roots are  {\em harmonically separated}, with cross ratio equal to $-1, 2$ or $1/2$.
Thus the corresponding points  can be transformed to vertices of a square. We can use
the remaining freedom in the M\"obius transformations to set
two roots to $\pm i$. Then either $b_1=0$ in which case
the remaining two roots are at $\pm 1$, or  $b_1=\pm 
1$ in which case
the roots are at $2\pm i$ or $-2\pm i$. Both squares are equivalent under the
Mobi\"us transformations. Thus the resulting five points form a square pyramid (Figure 3).
\koniec
%There are examples of non--maximally separated configurations
%which only allow one root in (\ref{quadratic}). This will occur
%either if $(a=0, b\neq 0)$, or if the discriminant $b^2-ac$ vanishes.
%Both cases correspond to a one--parameter family of configurations
%of four points on the plane. The former one is always equivalent
%to $\pm i, c, c+2$, where $c$ is an arbitrary 
%point.
%%%%%%%%%%%%%%%%%%%%%%%%%%%%%%%%%%%%%%%%%%%%%%%%%%%%%%%%%%%%%%%%%%%%%%%%%%%%%%
\section{Motivation: Twistor theory of $G_2$ structures}
We shall close the paper explaining the motivation of characterising binary sextics with 
$\mathcal{I}=0$ coming from twistor theory of $G_2$ structures. For the sake of the following discussion
twistor theory is a correspondence \cite{Pe76} between global algebraic geometry of curves in complex two--folds
or three-folds, and local differential geometry on the moduli spaces of these curves. 

Let ${\mathcal Z}$ be a complex two--fold with a family of rational curves $L_m\cong\CP^1$
parametrised by points $m\in M$, where $M$ is some complex manifold. For any $m$, the embedding of
$L_m$ in ${\mathcal Z}$ is, to the first order, described by the normal bundle
$N(L_m)=T{\mathcal Z}/TL_m$. This is a holomorphic line bundle ${\mathcal{O}}(k)$, for some integer $k$ which
we shall assume to be positive and even. The obstruction group $H^1(L_m, N(L_m))=0$ vanishes,
and the Kodaira deformation theorem \cite{kod}  states that there exists a canonical isomorphism
\[
T_m M\cong H^0(L_m, N(L_m))=\mbox{Sym}^k(\C^{2\vee})
\]
between vectors tangent to $M$ and binary quantics of degree $k$. This is where the quadratic invariant
(\ref{qin})  becomes relevant: if $k$ is even, we can  define a holomorphic conformal 
structure $[g]$ on $M$ by declaring a  vector field $U\in\Gamma (TM)$ to be null iff the 
corresponding quantic has $\mathcal{I}_2(U)=0$. This is a quadratic condition, and the holomorphic 
light-cone of any point in $M$ is a surface of co--dimension one in $M$, so $[g]$ is indeed well
defined.

Let us now restrict to the case $k=6$, where vector fields correspond to binary sextics, and 
$\mbox{dim}_\C M=7$, \cite{DG10} . In this case there exists a skew--symmetric 
three--form $\Psi\in \Lambda^3(M)$ given by
\[
\Psi(U, V, W)=<<U, V>_3, W>_3
\]
where $<, >_3$ is the third transvectant (\ref{k_th_trans})
of two binary sextics, and we use the same symbols
to denote vector fields and corresponding sextics. This three--form is compatible with $[g]$
is a sense that
\[
(U\hook \Psi)\wedge (U\hook\Psi)\wedge\Psi=0, \quad\mbox{iff}\quad {\mathcal{I}_2}(U)=0.
\]
The invariants ${\mathcal I}_2$ and $\Psi$ have weights six and nine respectively,
so changing a metric $g\in [g]$ yields $g\rightarrow \Omega^6 g, \Psi\rightarrow \Omega^9\Psi$
where $\Omega$ is a non--vanishing function on $M$. Therefore the structure group of $TM$ reduces
to $GL(2, \C)\subset {G_2}^{\C}\times \C$, where ${G_2}^{\C}$ is the complexification of the
 exceptional Lie group $G_2$.

There are only few known examples of this construction which lead to positive--definite 
$G_2$ structures on Riemannian manifolds $M_\R$. Bryant's weak $G_2$ holonomy metric \cite{Bryant1} on 
$M_\R=SO(5)/SO(3)$ arises from a family of $Sp(4)$ invariant rational sextics \cite{DS10}.
Another example corresponds to ${\mathcal Z}=\CP^2$, and $M$ being the homogeneous space
$PSL(3, \C)/\C^*$ of ternary cuspidal cubics in ${\mathcal Z}$. The cuspidal cubics are
rational, but singular. The Kodaira theory nevertheless  applies as the contact lifts
of the cuspidal cubics to $T\CP^2$ are smooth. There exist three real slices of $M$, one of which 
is $M_\R=SU(2, 1)/U(1)$. In \cite{DD} it is shown that $M_\R$ admits a $G_2$ structure which is co--calibrated
\cite{DD}.

\section{Conclusions}
We have found a geometric interpretation of vanishing of the quadratic
invariant associated with a binary sextic, and more generally with any binary quantic of even degree. The result should have an interpretation in the theory of hyperelliptic curves, as
projective equivalence classes of binary sextics with distinct roots correspond to points
on the moduli space of genus two algebraic curves. Igusa \cite{igusa} considered a zero locus of various invariants for sextics, but the case of vanishing quadratic invariant
has not been analysed in this work.
 
The problem addressed in this paper can also be reformulated
in the context of the representation theory of $SO(3, \C)$
which is related to $SL(2, \C)$ by the homomorphism $SL(2, \C)/\Z_2\cong SO(3, \C)$.
The isomorphism $\C^3=\mbox{Sym}^2(\C^{2\vee})$ identifies complex vectors
in $\C^3$ with symmetric $2$ by $2$ matrices with complex coefficients.
The {\em null} vectors correspond to rank one matrices, and this
gives rise to a holomorphic 
conformal structure on $\C^3$. The seven--dimensional space 
$\mbox{Sym}^6(\C^{2\vee})$ of binary sextics
is identified with a subspace of $\mbox{Sym}^3(\C^{3\vee})$ which consist
of harmonic ternary cubics, i. e. those forms $\Psi_{ijk}Z^iZ^jZ^k,
i, j, k=1, \dots, 3$ which satisfy $\delta^{ij}\Psi_{ijk}=0$.
Hitchin \cite{hitchin} showed that  a generic harmonic ternary 
cubic in $\CP^2$ passes through two sets of six points corresponding to
six axes of a  regular icosahedron in $\C^3$. In this formulation
the  quadratic invariant is given by the norm of the cubic 
$\Psi_{ijk}$ taken w. r. t. the conformal structure defined above.
It is however not clear what is the geometric meaning of its vanishing
in terms of Hitchin's icosahedron.
%\section*{Appendix}
\section*{Appendix: Two component spinors}
\appendix
\setcounter{equation}{0}
\def\theequation{\thesection{A}\arabic{equation}}
\label{sec_spinors}
A convenient way to represent binary quantics and the associated
invariants uses the two--component spinor notation \cite{PR}. 
Let the capital letters $A, B, \dots$ denote indices taking values
$0$ and $1$. The general quantic  (\ref{quantic1}) is represented by a 
{\it symmetric spinor} of valence $2n$. The Fundamental 
Theorem of Algebra states that any such spinor factorises
into valence one spinors
\[
\psi_{AB\cdots C}=\alpha_{(A}\beta_B\dots\gamma_{C)},
\]
where the round brackets on the RHS denote the symmetrisation.
The binary  quantic (\ref{quantic1}) is then given by
\begin{eqnarray*}
\psi&=&\psi_{AB\dots C}\pi^A\pi^B\dots \pi^C\\
&=&(\alpha_0 x+\alpha_1 y)(\beta_0 x+\beta_1 y)\dots
(\gamma_0 x+\gamma_1 y)\sim (x-x_1)(x-x_2)\dots (x-x_{2n}),
\end{eqnarray*}
where $\pi^0=x, \pi^1=y$ and
$
\psi_0=\psi_{00\dots 0}, \quad\psi_1=\psi_{10\dots 0},\quad \dots, \quad
\psi_{2n}=\psi_{11\dots 1}.
$
Thus the complex numbers $x_1=-\alpha_1/\alpha_0, x_2=-\beta_1/\beta_0, \dots, x_{2n}=-\gamma_1/\gamma_0$ are the roots of
the inhomogeneous polynomial of degree $2n$ obtained by setting $y=1$.
The invariant (\ref{qin}) is in this notation given by
\[
{\mathcal I}=\psi_{AB\dots C}\psi^{AB\dots C},
\]
where the indices are lowered by the anti--symmetric matrix 
$\varepsilon_{AB}$ with $\varepsilon_{01}=1$, so that
\[
\psi_{AB\dots C}=\psi^{PQ\dots R}\varepsilon_{PA}\varepsilon_{QB}\dots
\varepsilon_{RC}.
\]
The $k$th transvectant (\ref{k_th_trans}) is
\[
<\psi, \phi>_k=\varepsilon^{A_1B_1}\dots \varepsilon^{A_kB_k}
\psi_{A_1\dots A_kA_{k+1}\dots A_m}\phi_{B_1\dots B_kB_{k+1}\dots B_n}\pi^{A_{k+1}}
\dots \pi^{A_m}\pi^{B_{k+1}}\dots\pi^{B_n}.
\]
Using the spinor notation gives a simple proof of the following algebraic interpretation
of the condition ${\mathcal I}=0$
\begin{lemma} An even degree quantic $\psi\in V_{2n}$ with 
distinct roots is a sum 
of ${(2n)}^{th}$ powers of its factors iff
$\mathcal{I}(\psi)=0$.
\end{lemma}\noindent
{\bf Proof.}
%The condition stated in the Lemma is 
%\be
%\label{elliot0}
%\psi=c_1(x-x_1)^{2n}+c_2(x-x_2)^{2n}+\dots+c_{2n}(x-x_{2n})^{2n}
%\ee
Let $\alpha_A, \beta_A, \dots, \gamma_A$ be homogeneous 
coordinates of the points in $\CP^1$ corresponding to the roots of $\psi$
so that the condition stated in the Lemma 
becomes\footnote{Formula (\ref{elliot}) holds for any binary quantic of odd degree, as then ${\mathcal{I}}$ vanishes identically.}
\be
\label{elliot}
\alpha_{(A}\beta_B\dots \gamma_{C)}=c_1\;\alpha_{A}\alpha_B\dots \alpha_C
+c_2\;\beta_{A}\beta_B\dots \beta_C+\dots 
+c_{2n}\;\gamma_{A}\gamma_B\dots \gamma_C,
\ee
where $c_1, c_2, \dots, c_{2n}$ are some constants determined by 
the roots. The condition $\mathcal{I}=0$ is then equivalent to
\[
\frac{{\p^{2n}}\psi}{\p\alpha \p\beta\dots \p\gamma}=0, \quad
\mbox{where}\quad \frac{\p}{\p \alpha}:=\alpha^A\frac{\p}{\p \pi^A}\;\mbox{etc.}
\]
Thus 
\[
\frac{{\p^{2n-1}}\psi}{\p\beta\dots \p\gamma}=\tilde{c}_1 (\alpha_A \pi^A)
\]
for some constant $\tilde{c}_1$, as both sides are 
homogeneous of degree one in
$\pi^A$.  The successive $(2n-2)$ integrations yield (\ref{elliot}), with
$
c_1={(2n)!}^{-1}{<\alpha, \beta>_1}^{-1}
\dots {<\alpha, \gamma>_1}^{-1}(\tilde{c}_1)$
 etc.
\koniec

\end{document}